\title{\vspace{-0.6cm} A counterexample to the Alon-Saks-Seymour conjecture\\ and related problems}
\author{Hao Huang \thanks{Department of Mathematics, UCLA, Los
Angeles, CA, 90095. Email: huanghao@math.ucla.edu.}
\and
Benny Sudakov\thanks{Department of Mathematics,
UCLA,  Los Angeles, CA 90095. Email: {\tt bsudakov@math.ucla.edu}.
Research supported in part by NSF CAREER award DMS-0812005 and by a USA-Israeli BSF grant.}}

\documentclass[11pt]{article}
\usepackage{amsfonts,amssymb,amsmath,latexsym}
\usepackage{verbatim}

\date{}

\newtheorem{thm}{Theorem}[section]
\newtheorem{prop}[thm]{Proposition}
\newtheorem{lemma}[thm]{Lemma}
\newtheorem{cor}[thm]{Corollary}
\newtheorem{claim}[thm]{Claim}
\newtheorem{conj}[thm]{Conjecture}
\newtheorem{ques}[thm]{Question}

\newenvironment{pf}
      {\medskip\noindent{\bf Proof.}\hspace{1mm}}
      {\hfill$\Box$\medskip}

\def\qed{\ifvmode\mbox{ }\else\unskip\fi\hskip 1em plus 10fill$\Box$}

\newcommand{\bp}{\textup{\textbf{bp}}}
\begin{document}
\maketitle

\begin{abstract}
Consider a graph obtained by taking edge disjoint union of $k$ complete bipartite graphs.
Alon, Saks and Seymour conjectured that such graph has chromatic number at most $k+1$.
This well known conjecture remained open for almost twenty years.
In this paper, we construct a counterexample to this
conjecture and discuss several related problems in combinatorial geometry
and communication complexity.
\end{abstract}

\section{Introduction}
\label{section_introduction}

Tools from linear algebra have many striking applications in the study of combinatorial problems.
One of the earliest such examples is the theorem of  Graham and Pollak \cite{MR0332576}.
Motivated by a communication problem that arose in connection with data transmission, they
proved that the edge set of a complete graph $K_k$ cannot be
partitioned into disjoint union of less than $k-1$ complete
bipartite graphs. Their original proof used Sylvester law of inertia. Over the years, this elegant result attracted a lot of attention and
by now it has several different algebraic proofs, see \cite{babai1992linear,
MR743808, MR679606, MR2407919}. On the other hand no purely combinatorial proof of this statement is known.

A natural generalization of Graham-Pollak theorem is to ask whether the same estimate holds
also for all graphs with chromatic number $k$. This problem was raised twenty years ago by
Alon, Saks and Seymour who made the following conjecture (see, e.g., survey of J.~Kahn, \cite{MR1319170}).

\begin{conj}
\label{ass}
If the edges of a graph $G$ can be partitioned into $k$ edge disjoint complete bipartite graphs, then
the chromatic number of $G$ is at most $k+1$.
\end{conj}
This question is also related to another long-standing open problem by
Erd\H{o}s, Faber and Lov\'{a}sz. They conjectured that
the edge disjoint union of $k$ complete graphs of order $k$ is $k$-chromatic.
Indeed, by replacing cliques in this problem by complete bipartite graphs
we obtain the Alon-Saks-Seymour conjecture. The question of Erd\H{o}s, Faber and Lov\'{a}sz  is still open. On the other hand, Kahn \cite{MR1141320} proved the 
asymptotic version
of their conjecture, showing that the chromatic number of  edge disjoint union of $k$ complete graphs of order $k$ has chromatic number at most $(1+o(1))k$.

Let $\bp(G)$ be the minimum number of bicliques (i.e., complete
bipartite graphs) needed to partition the edges of graph $G$ and
$\chi(G)$ be the chromatic number of $G$. The Alon-Saks-Seymour
Conjecture can be restated as $\bp(G) \geq \chi(G)-1$. Until
recently, there was not much known about this conjecture. Using 
folklore result that the chromatic number of the union of
graphs is at most the product of their chromatic numbers, one can
easily get a lower bound $\bp(G) \geq \log_2 \chi(G)$. In
\cite{mubayi2009biclique}, Mubayi and Vishwanathan improved the
lower bound to $2^{\sqrt{2 \log_2 \chi(G)}}$. This estimate can be
also deduced from the well known result of Yannakakis
\cite{yannakakis1991expressing} in communication complexity. This
connection to communication complexity was discovered by Alon and
Haviv \cite{alonhaviv} (see Section 4 for details). Gao, McKay, Naserasr and Stevens
\cite{naserasr} introduced a reformulation of the Alon-Saks-Seymour
conjecture and verified it for graphs with chromatic number $k \leq
9$. The main aim of this paper is to obtain a superlinear gap
between chromatic number and biclique partition number, which
disproves the Alon-Saks-Seymour conjecture.

\begin{thm}
\label{main}
There exists an infinite collection of graphs $G$ such that
$\chi(G) \geq c\,\big(\bp(G)\big)^{6/5}$, for some fixed constant $c>0$.
\end{thm}

The study of (two-party) communication complexity, introduced by Yao \cite{yao1979some}, is an important topic in theoretical computer science
which has many applications.
In the basic model we have two players Alice and Bob who are trying to evaluate a
boolean function $f: X \times Y \rightarrow \{0, 1\}$. Alice only knows $x$, Bob only knows $y$ and thus they want to communicate with each other according to some fixed
protocol in order to compute $f(x,y)$. The goal is to minimize the amount of communication during the protocol.
The deterministic communication complexity $D(f)$ is the number of bits that needs to be exchanged for the worst inputs $x,y$  by the best protocol for $f$.
Let $M$ be a matrix of $f$, i.e., $M_{x,y}=f(x,y)$ and let $rk(M)$ be the rank of $M$.
It's known that $D(f) \geq \log_2 rk(M)$.
Lov\'{a}sz and Saks \cite{lovasz1993lattices} conjectured that this bound is not very far from
being tight. More precisely, their  {\em log-rank conjecture} says that
$D(f) \leq (\log_2 rk(M)^{O(1)}$. This problem is directly related to
the {\em rank-coloring conjecture} of Van Nuffelen \cite{van38rank} and Fajtlowicz
\cite{fajtlowicz1987conjectures} in graph theory. This conjecture, which was disproved by Alon and Seymour \cite{alon1989counterexample}, asked whether the chromatic number
of a graph $G$ is bounded by the rank of its adjacency matrix $A_G$.
It is known that separation result between $D(f)$ and $\log_2 rk(M)$ give corresponding separation between
$\chi(G)$ and $rk(A_G)$. Several authors gave such separation results, e.g.,  \cite{MR1189860, raz1995log}.
So far, the largest gap was obtained by Nisan and Wigderson \cite{nisan1995rank} who constructed an
infinite family of matrices such that $D(f)>(\log_2 rk(M))^{\log_2 3}$.

Similar to the rank-coloring problem, the Alon-Saks-Seymour conjecture is also
closely related to a well known open problem in communication complexity. This communication problem, which is called
{\em clique versus independent set} ($CL\text{-}IS$ for brevity), was introduced by Yannakakis \cite{yannakakis1991expressing} in 1988.
In this problem, there is a publicly known graph $G$, Alice gets a clique $C$ of $G$ and Bob gets
an independent set $I$ of $G$. Their goal is to output $|C\cap I|$, which is clearly either $0$ or $1$.
We will discuss connection between this problem and the Alon-Saks-Seymour conjecture and  show that our counterexample yields a first nontrivial
lower bound on the non-deterministic communication complexity of $CL\text{-}IS$  problem.

The rest of this short paper is organized as follows. In the next section
we describe a counterexample to the Alon-Saks-Seymour
Conjecture. In Section \ref{section_neighborlyfamily}, we consider minimal coverings of a graph by
bicliques in which every edge of the graph is covered at least once and at most $t$ times, for some parameter $t$. This
more general notion is closely related to the question in combinatorial geometry about a neighborly family of boxes. We show that a natural variant of
the Alon-Saks-Seymour conjecture for this more general parameter fails as well.
In Section \ref{section_communicationcomplexity}, we discuss connections with
communication complexity and use our counterexample to
obtain a new lower bound on nondeterministic communication complexity of
clique vs. independent set problem. The final section contains some
concluding remarks and open problems.

\vspace{0.25cm}
\noindent
{\bf Notation.}\, The $n$-dimensional cube $Q_n$ is $\{0,1\}^n$ and
two vertices $x, y$ of $Q_n$ are adjacent $x \sim y$ if and only
if they differ in exactly one coordinate. A $k$-dimensional subcube of $Q_n$
is a subset of $\{0,1\}^n$ which can be written as $\{x=(x_1,\cdots,x_n) \in
Q_n:x_i=a_i,~\forall i \in T\}$, where $T$ is a set of $n-k$ coordinates
(called fixed coordinates), each $a_i$ is a fixed element in $\{0,1\}$. In
addition, we write $1^n$ and $0^n$ to represent the all-one and all-zero
vector in $Q_n$ and use $Q_n^{-}$ to indicate the set $Q_n \backslash
\{1^n,0^n\}$. Given two subset $X \subset Q_k$ and $Y \subset Q_{\ell}$ we denote by
$X\times Y$ a subset of cube $Q_{k+\ell}$ which consists of all binary vectors
$(x, y)$ with $x \in X$ and $y\in Y$.

For graph $G=(V,E)$ with vertex set $V$ and edge set
$E$, we denote by $\chi(G),~ \alpha(G),~ \bp(G)$ the chromatic number,
independence number and biclique partition number respectively. The collection of all independent sets in $G$ is denoted
by $\mathcal{I}(G)$. Similarly $\mathcal{C}(G)$ stands for the set of all cliques in $G$. The $OR$ product of two graphs $G$ and
$H$ is defined as a graph with vertex set equal to the Cartesian
product $V(G) \times V(H)$, two vertices $(g,h)\sim (g',h')$ iff $g
\sim g'$ in $G$ or $h \sim h'$ in $H$. The $m$-blowup of a graph $G$
is obtained by replacing every vertex $v$ of $G$ with an independent set $I_v$ of size $m$
and by replacing every edge $(u,v)$ of $G$ with a complete bipartite graph, whose
parts are the independent sets $I_u$ and $I_v$.
We also use the notation $\mathcal{B}(U,W)$ to indicate a biclique with two
parts $U$ and $W$.

To state asymptotic results, we utilize the following standard
notations. For two functions $f(n)$ and $g(n)$, write $f(n)=
\Omega(g(n))$ if there exists a positive constant $c$ such that
$\lim \inf_{n \rightarrow \infty} f(n)/g(n) \geq c$, $f(n) =
o(g(n))$ if $\lim \sup_{n \rightarrow \infty} f(n)/g(n) = 0$. Also,
$f(n) = O(g(n))$ if there exists a positive constant $C > 0$ such
that $\lim \sup_{n \rightarrow \infty} f(n)/g(n) \leq C$.

\section{Main Result}
\label{section_mainresult}
In this section we describe a counterexample to the Alon-Saks-Seymour conjecture.
Our construction is inspired by and is somewhat similar to Razborov's counterexample
to the rank-coloring conjecture \cite{MR1189860}.
Consider the following graph $G=(V,E)$. Its vertex set is
$V(G)=[n]^7=\{(x_1, \cdots,
x_7): x_i \in [n]\}$. For any two vertices $x=(x_1, \cdots, x_7)$,
$y=(y_1, \cdots, y_7)$ in $V(G)$, let $\rho$ be the comparing function which records all
coordinates in which they differ. More precisely,
$\rho(x,y)=(\rho_1(x,y), \cdots, \rho_7(x,y)) \in Q_7$, such that
\begin{equation*}
\rho_i(x,y)=\begin{cases}
1 & \textrm{if}~x_i \neq y_i \\
0 & \textrm{if}~x_i = y_i
\end{cases}
\end{equation*}
Two vertices $x$ and $y$ are adjacent in $G$  if and only if
$\rho(x,y) \in S$, where  $S$ is the following subset of the cube $Q_7$
$$ S= Q_7 \, \backslash\, \big[(1^4 \times Q_3^{-}) \cup \{0^4 \times
0^3\} \cup \{0^4 \times 1^3\}\big].$$
In the rest of this section we show that this graph $G$ satisfies the assertion of Theorem \ref{main}.

\begin{prop}\label{independentnumber}
The independence number of $G$ satisfies $\alpha(G)=O(n)$.
\end{prop}

\begin{pf}
Let $I$ be an independent set in $G$. For any set of indices $T=\{i_1, \ldots, i_t\}
\subset \{1,2, \cdots,7\}$, let $p_T$ be the natural projection of $[n]^7$ to $[n]^T$.
For every vector $x \in [n]^7$ it outputs the restriction of $x$ to the coordinates in $T$, i.e.,
$p_T(x)=(x_{i_1}, \ldots, x_{i_t})$.
For convenience, we will for example write $p_{1234}$ instead of
$p_{\{1,2,3,4\}}$. It is easy to check from the definition of $S$, that any two vertices $x,y \in G$ which agree on one of the first 4 coordinates and satisfy
$p_{1234}(x)\not = p_{1234}(y)$ are adjacent in $G$. Hence, any two vectors in
$p_{1234}(I)$ differ in all their coordinates and therefore $|p_{1234}(I)|
\leq n$. If in addition, we also have  for every element $x \in p_{1234}(I)$,
$|p_{1234}^{-1}(x) \cap I| \leq 3$, then $|I| \leq 3|p_{1234}(I)|
=O(n)$ and the proof is complete.

Otherwise, we may assume the existence of $\widetilde{x} \in [n]^4$
and different vertices $\widetilde{x}_1, \widetilde{x}_2,
\widetilde{x}_3, \widetilde{x}_4 \in I$ such that
$p_{1234}(\widetilde{x}_i)=\widetilde{x}$. By the definition of $S$,
it is easy to see that
$p_{567}(\widetilde{x}_i)$ differ in every coordinate.
Since $1^7 \in S$, we have that any two vertices of $G$ which differ in all $7$ coordinates are
adjacent. This implies that if there is a vertex $z \in I$ with $p_{1234}(z)$ different from
$\widetilde{x}$, then $p_{567}(z)$ and $p_{567}(\widetilde{x}_i)$
are equal in at least one coordinate. Since the number of coordinates of $p_{567}(I)$ is only $3$ and there are $4$ vertices $\widetilde{x}_1, \widetilde{x}_2,
\widetilde{x}_3, \widetilde{x}_4$, we have that two of these vertices agree with $p_{567}(I)$ (and hence with each other) in the same coordinate.
This contradicts the fact that $p_{567}(\widetilde{x}_i)$ differ in all coordinates, and implies that there is
only one element in $p_{1234}(I)$. Again, by the definition of $S$, the vertices in $I$ are different in each of the last three
coordinates. As a result $|I|=|p_{567}(I)| \leq n$.
\end{pf}

\begin{cor}\label{chromaticnumber}
The chromatic number of $G$ is at most $\Omega(n^6).$
\end{cor}
\begin{pf}
Apply Proposition \ref{independentnumber} together with  the well-known fact that $\chi(G) \geq
\frac{|V(G)|}{\alpha(G)}$.
\end{pf}

\begin{prop} \label{bpnumber}
The biclique partition number satisfies $\bp(G)=O(n^5)$.
\end{prop}
Before going into the details of the proof of this statement, first we need the
following two lemmas.
\begin{lemma} \label{lemmadecomposition}
S can be partitioned into disjoint union $S=\cup_{i=1}^{30}S_i$, where
$S_i$ are $2$-dimensional subcubes in $Q_7$.
\end{lemma}
\begin{pf}
We need the following simple observations.

\noindent $(a)$~$Q_3^{-}$ is a disjoint union of
$1$-dimensional subcubes.\\
$(b)$~$Q_3$ can be decomposed into disjoint union of $2$-dimensional subcubes.\\
$(c)$~For every $R_1 \subset Q_4$,  the set $R_1 \times Q_3$
can be decomposed into disjoint union of $2$-dimensional subcubes.\\
$(d)$~For any $x_1 \sim x_2$ in $Q_4$, $y_1 \sim y_2$ in $Q_3$,
$(x_1,y_1),(x_1,y_2),(x_2,y_1),(x_2,y_2)$ is a
$2$-dimensional subcube in $Q_7$.\\
$(e)$~For any $x_1 \sim x_2$ in $Q_4$, $(x_1 \times Q_3^{-}) \cup
(x_2 \times Q_3^{-})$ can be decomposed into disjoint union of
$2$-dimensional subcubes.

To verify $(a)$ note that, $Q_3^{-}=\{(0,0,1),(0,1,1)\} \cup \{(0,1,0),(1,1,0)\}
\cup \{(1,0,0),(1,0,1)\}$. Claims $(b)$ and $(d)$ are obvious by the definition
of cubes. Claims $(c)$ is an immediate corollary of $(b)$, and claim $(e)$ follows easily from
$(a)$ and $(d)$.

Next we can partition the set $S=Q_7\, \backslash\, \big[(1^4 \times Q_3^{-})
\cup \{0^4 \times 0^3\} \cup \{0^4 \times 1^3\}\big]$ into the following
3 disjoint subsets $S', S'', S'''$ and show that each of them is itself a disjoint union of $2$-dimensional subcubes.

\begin{equation*}
S'=\begin{cases}
(0,0,0,0)\times Q_3^{-}\, \cup \,(0,0,0,1) \times Q_3^{-}\\
(0,0,1,1)\times Q_3^{-}\, \cup\, (1,0,1,1) \times Q_3^{-}\\
(0,1,0,1)\times Q_3^{-}\, \cup \,(0,1,1,1) \times Q_3^{-}\\
(1,1,0,1)\times Q_3^{-} \,\cup \,(1,0,0,1) \times Q_3^{-}
\end{cases}
\end{equation*}
This set can be partitioned into disjoint union of $2$-dimensional subcubes, using claim $(e)$.
\begin{equation*}
S''=\begin{cases}
(1,1,1,1)\times 0^3 \, \cup \, (1,1,0,1) \times 0^3 \, \cup \, (1,0,1,1) \times 0^3 \, \cup \, (1,0,0,1) \times 0^3\\
(1,1,1,1)\times 1^3 \, \cup \, (1,1,0,1) \times 1^3 \, \cup \,(1,0,1,1) \times 1^3 \, \cup \,(1,0,0,1) \times 1^3\\
(0,1,1,1)\times 0^3 \, \cup \, (0,1,0,1) \times 0^3 \, \cup \,(0,0,1,1) \times 0^3 \, \cup \, (0,0,0,1) \times 0^3\\
(0,1,1,1)\times 1^3 \, \cup \, (0,1,0,1) \times 1^3 \, \cup \, (0,0,1,1) \times 1^3 \, \cup \,(0,0,0,1) \times 1^3
\end{cases}
\end{equation*}
Note  that, every line in the definition of $S''$ describes a $2$-dimensional subcube. This shows that $S''$ is
a disjoint union of four $2$-dimensional subcubes.

\begin{equation*}
S'''= \begin{cases} (0,0,1,0)\times Q_3 \, \cup \, (0,1,0,0)\times Q_3 \, \cup \, (1,0,0,0)\times Q_3 \, \cup \, (0,1,1,0)\times Q_3\\
(1,0,1,0)\times Q_3 \, \cup \, (1,1,0,0)\times Q_3 \, \cup \, (1,1,1,0)\times Q_3
\end{cases}
\end{equation*}
To decompose this set into disjoint union of $2$-dimensional subcubes, one can use claim $(c)$.

Finally, it is easy to verify that indeed $S=S' \cup S'' \cup S'''$ and hence
$S$ can be partitioned into $2$-dimensional subcubes.
\end{pf}

Using the decomposition $S=\cup_{i=1}^{30}S_i$ from Lemma
\ref{lemmadecomposition}, we can define the following subgraphs $G_i \subset G$.
The vertex set $V(G_i)=V(G)$ and two vertices $x,y \in G_i$ are
adjacent if and only if $\rho(x,y) \in S_i$. From this definition,
it is easy to see that $G$ is the edge disjoint union of
subgraphs $G_i$. Next we will show that every $G_i$ has a small
biclique partition number.

\begin{lemma}\label{lemmabp}
$\bp(G_i) \leq n^5$.
\end{lemma}
\begin{pf}
Recall that the set $S_i$, which is used to define edges of $G_i$, is a $2$-dimensional
subcube of $Q_7$. Therefore there exists a set $T= \{t_1, \ldots, t_5\} \subset \{1,
\cdots , 7\}$ of fixed coordinates and $a_1, \ldots, a_5 \in \{0,1\}$, such that
$S_i=\big\{x=(x_1, \cdots, x_7): x_{t_j}=a_j,~\forall\, 1 \leq j \leq 5\big\}$. Now we
define graph $\widetilde{G}_i$. Its vertex set $V(\widetilde{G}_i)=[n]^5$ and two vertices
$\widetilde{x}$ and $\widetilde{y}$ are adjacent in
$\widetilde{G}_i$ if an only if $\rho(\widetilde{x},\widetilde{y})=(a_1, \ldots, a_5)$. It is rather straightforward
to see that $G_i$ is a $n^2$-blowup of $\widetilde{G}_i$.

To complete the proof of this lemma we need two basic facts
about biclique partition number. The first one says that
for any graph $H$, $\bp(H) \leq |V(H)|-1$. Indeed, removing stars rooted
at every vertex, one by one, we can partition every graph on $h$ vertices into
$h-1$ bicliques. The second one, claims that if $H$ is a
blowup of $\widetilde{H}$, then $\bp(H) \leq \bp(\widetilde{H})$.
To prove this, note that the blowup of biclique is a biclique itself. Therefore blowup of
all the bicliques in a partition of $\widetilde{H}$
becomes a biclique partition of $H$.

These two statements, together with the fact (mentioned above) that $G_i$ is the blowup of $\widetilde{G}_i$,
imply that $\bp(G_i) \leq \bp(\widetilde{G}_i) \leq |V(\widetilde{G}_i)|-1 \leq n^5$.
\end{pf}

\vspace{0.2cm}
\noindent
{\bf Proof of Proposition \ref{bpnumber}.}\hspace{1mm}
Using that $G$ is the edge disjoint union of $G_i$ together with Lemma \ref{lemmabp}, we
conclude that, $\bp(G)=\bp(\cup_{i=1}^{30} G_i) \leq
\sum_{i=1}^{30} \bp(G_i) = O(n^5)$. \hfill$\Box$\medskip

Propositions \ref{chromaticnumber} and \ref{bpnumber} show that graph $G$,
which we constructed, indeed satisfies the assertion of Theorem \ref{main}
and disproves the Alon-Saks-Seymour Conjecture.

\section{Neighborly family of boxes and  $t$-biclique covering number}
\label{section_neighborlyfamily}
The Alon-Saks-Seymour conjecture deals with the minimum number of bicliques
needed to cover all the edges of a given graph $G$ exactly once.
It is also  very natural to consider a more general
problem in which we are allowed to cover the edges of graph
at most $t$ times. A {\em $t$-biclique covering} of a graph $G$ is a collection of
bicliques that cover every edge of $G$ at least once and at most $t$ times. The minimum size of such covering is called
the {\em $t$-biclique covering number} and is denoted by  $\bp_t(G)$. In particular,
$\bp_1(G)$ is the usual biclique partition number $\bp(G)$.

In addition to being an interesting parameter to study in its own right, the $t$-biclique covering number is also closely
related to the question in combinatorial geometry about neighborly family of boxes. A finite family
$\mathcal{C}$ of $d$-dimensional convex polytopes is
called {\em $t$-neighborly} if $d-t \leq dim(C \cap C') \leq d-1$ for
every two distinct members $C$ and $C'$ of $\mathcal{C}$. One
particularly interesting case is when $\mathcal{C}$ consists of
$d$-dimensional boxes with edges parallel to the coordinate axes.
This type of box is  called {\em standard box}. Using
Graham-Pollak theorem, Zaks \cite{zaks1979bounds} proved  that the
maximum possible cardinality of a $1$-neighborly family of standard
boxes in $\mathbb{R}^d$ is precisely $d+1$.
His result was generalized  by Alon
\cite{alon1997neighborly}, who proved that $\mathbb{R}^d$ has a $t$-neighborly family of $k$
standard boxes if and only if the complete graph $K_k$ can has $t$-biclique covering of size $d$.
This shows that the problem of determining the maximum possible cardinality of $t$-neighborly families of standard boxes and the problem of
computing the $t$-biclique covering number of complete graphs are equivalent.

In his paper \cite{alon1997neighborly}, Alon gave asymptotic estimates for $\bp_t(K_k)$ showing that
$$(1+o(1))\big(t!/2^t\big)^{1/t}k^{1/t}\leq \bp_t(K_k) \leq (1+o(1))t\,k^{1/t}.$$
There is still gap between these two bounds and the problem of determining the right constant before $k^{1/t}$ is wide open even for
the case when $t=2$. Using a different proof, we obtain here a slightly better lower bound of order roughly  $\big(t!/2^{t-1}\big)^{1/t}k^{1/t}$. 
For $t=2$ it improves the above estimate by a factor of $\sqrt{2}$.

\begin{prop}
If there exists a $t$-biclique covering of $K_k$ of size $d$, then
$k \leq 1+\sum_{s=1}^t 2^{s-1}{d \choose s}$.
\end{prop}
\begin{pf}
Suppose that the edges of $K_k$ are covered by the bicliques
$\{\mathcal{B}(U_j,W_j)\}_{j=1}^d$, such that every edge is covered at least once and at most $t$-times.
For every nonempty subset of indices $S \subset [d]$ of size $|S| \leq t$ let $H_S=\cap_{j \in S} \mathcal{B}(U_j,W_j)$ and let $A_S$ be the adjacency matrix of $H_S$.
Let $J$ be $k\times k$ matrix of ones and let $I$ be the $k\times k$ identity matrix. Then $J-I$ is the adjacency matrix of $K_k$ and it is easy to see, using
the inclusion-exclusion principle, that
$$J-I=\sum_{S \subset [d], 0<|S|\leq t} (-1)^{|S|} A_S.$$
Also note that for $|S|=s$, the graph $H_S$ is the disjoint union of at most $2^{s-1}$ smaller bicliques. Indeed, for every binary vector $z=(z_1, \ldots, z_{s-1})$
consider a complete bipartite graph with parts
$$X_z=\cap_{j, z_j=0}\,U_j\, \cap_{j, z_j=1}\,W_j \cap\ U_s ~~\mbox{and}~~Y_z=\cap_{j, z_j=0}\, W_j\, \cap_{j, z_j=1}\,U_j \cap W_s.$$
It is not difficult to check that these bicliques are disjoint and their union is $H_S$.
Therefore, for every $S \subset [d], 0<|S|=s\leq t$ we can write $A_S=\sum_i B_{i,S}$ where $B_{i,S}$ is an adjacency matrix of a biclique and $ 1 \leq i \leq 2^{s-1}$.
Thus  we obtain that $J-I$ can be written as a linear combination of at most $m=\sum_{s=1}^t 2^{s-1}{d \choose s}$ adjacency matrices of complete bipartite graphs.

Now to complete the proof we use the elegant trick of Peck \cite{MR743808} (we can use here other known proofs of Graham-Pollak theorem as well).
For bipartite graph with adjacency matrix $B_{i,S}$ let
$B'_{i,S}$ be $k\times k$ matrix which contains only ones in positions whose row index lies in the first part of the bipartition and whose column index lies in the second part of
the bipartition, the rest of the entries of $B'_{i,S}$ are zeros. Since the corresponding bipartite graph is complete, $B'_{i,S}$ has rank one. Furthermore, the matrix
$B_{i,S}-2B'_{i,S}$ is antisymmetric. As a result we can write $J-I$ as a linear combination of at most $m$ rank one matrices plus
some antisymmetric matrix $T$. Since an antisymmetric real matrix has only imaginary eigenvalues, $I+T$ must have a full rank $k$. But its rank can not exceed the
rank of the linear combination of at most $m$ rank one matrices plus $J$. As $J$ has rank one as well, this implies that
$k \leq m+1=1+\sum_{s=1}^t 2^{s-1}{d \choose s}$ and completes the proof.
\end{pf}

As we already mentioned in the introduction, the motivation for the Alon-Saks-Seymour conjecture
comes from Graham-Pollak theorem which says that $\bp(K_k) \geq k-1$. Similarly, based on the lower bound of Alon
that $\bp_t(K_k) \geq \Omega\big(k^{1/t}\big)$, one can consider the following very natural
generalization of this conjecture.

\begin{ques}
\label{gass}
Is it true that for every fixed integer $t>0$, there exist a constant $c=c(t)$ such
that $\bp_t(G) \geq c \big(\chi(G)\big)^{1/t}$ for all graphs $G$?
\end{ques}

Recall that in Section 2 we constructed a graph $G$ with $|V(G)|=n^7$ vertices such that $\alpha(G)=O(n)$ and $\bp(G)=O(n^5)$.
Consider the $OR$ product (defined in the introduction) of $t$ copies of $G$. We show that the graph $G^t$ gives a negative answer to the above question
for all positive integers $t$. This follows from the following sequence of claims.

\begin{claim} \label{tindnumber}
$\alpha(G^t) \leq \alpha(G)^t = O(n^t)$.
\end{claim}

\begin{pf}
We only need to prove $\alpha(G \times H) \leq \alpha(G) \alpha(H)$
for any two graphs $G$ and $H$, since  then the claim follows by induction on $t$. To
prove this statement, consider the maximum independent set $I \in G \times H$. Let
$I'=\{v \in G ~|~(v,u) \in I ~\mbox{for some}~u\in H\}$ be the projection of $I$ on $V(G)$.
By the definition of $OR$ product, this is an independent set in $G$ and therefore has size at most
$\alpha(G)$. Similarly, if $I''$ is the  projection of $I$ on $V(H)$ then $|I''| \leq \alpha(H)$.
To complete the proof note that $I$ is a subset of $I' \times I''$ and therefore its size cannot exceed $\alpha(G)\alpha(H)$.
\end{pf}

\begin{cor}\label{tchromaticnumber}
$\chi(G^t) = \Omega (n^{6t})$.
\end{cor}
\begin{pf}
By Claim \ref{tindnumber}, $\chi(G^t) \geq
\frac{|V(G^t)|}{\alpha(G^t)} \geq \frac{n^{7t}}{\alpha(G)^t} =
\Omega (n^{6t})$.
\end{pf}

\begin{claim} \label{tbpnumber}
$\bp_t(G^t) \leq t\,\bp(G)$.
\end{claim}
\begin{pf}
Consider graphs $H_i, 1 \leq i \leq t$ with vertex set $V(H_i)=V(G^t)$ such that two vertices
$(h_1,\cdots,h_t)$ and $(h'_1,\cdots,h'_t)$ are adjacent in $H_i$ if and only if
$h_i \sim h'_i$ in $G$. Note that $H_i$ is an $n^{t-1}$-blowup of $G$ and therefore
$\bp(H_i)=\bp(G)$. Also it is easy to see that every edge in $G^t$ is covered by some $H_i$.
Since the number of graphs $H_i$ is $t$, every edge of
$G^t$ is covered at most $t$ times. Then the union of minimum biclique partitions of all $H_i$
gives a $t$-biclique covering of $G$. Hence $\bp_t(G^t) \leq \sum_{i=1}^t \bp(H_i) \leq t\,\bp(G)$.
\end{pf}

\begin{claim}
$\bp_t(G^t) \leq c\big(\chi(G^t)\big)^{\frac{5}{6t}}$ for some constant $c=c(t)$.
\end{claim}
\begin{pf}
By Claims \ref{tchromaticnumber} and \ref{tbpnumber}, $\bp_t(G^t) \leq
t\, \bp(G) = O(tn^5) \leq c(t)\big (\chi(G^t)\big)^{\frac{5}{6t}}$.
\end{pf}

\noindent
This shows that the answer to the Question \ref{gass} is negative for all natural $t$.

\section{Clique vs. independent set communication problem}
\label{section_communicationcomplexity} 
In the introduction, we already defined the two-party communication model
and discussed the concept of deterministic communication complexity. 
Here we need a few additional notions and definitions
(see e.g., \cite{MR1426129} for more details).  
The {\em non-deterministic communication complexity} $N^1(f)$ of a function $f$ is 
the smallest number of bits needed by an all powerful prover to convince 
Alice and Bob that $f(x,y)=1$. It is known that $N^1(f)=\lceil \log_2 C^1(f)\rceil$, where
$C^{1}(f)$ is the minimum number of monochromatic combinatorial rectangles needed to cover the
$1$-inputs of communication matrix $M$ of $f$ (recall that $M_{x,y}=f(x,y)$). With slight abuse of notation we will later
write $C^{1}(M)$ instead of $C^{1}(f)$. The numbers $N^0(f)$, $C^{0}(f), C^0(M)$ are defined similarly, 
and the relation $N^0(f)=\lceil \log_2 C^0(f)\rceil$ holds as well.

In this section we consider the communication complexity of the clique versus independent set problem
($CL\text{-}IS$). In this problem, there is a publicly known graph $\Gamma$, Alice gets a clique $C$ of $\Gamma$ and Bob gets
an independent set $I$ of $\Gamma$. Their goal is to output $|C\cap I|$, which is clearly either $0$ or $1$. This problem was first introduced by 
Yannakakis \cite{yannakakis1991expressing}, who also proposed the following algorithm to solve it. 
Given a graph $\Gamma$ on $m$ vertices, Alice sends to Bob a name of the vertex $v$
in $C$ whose degree in $\Gamma$ is at most $m/2$. Note that in this case we can reduce the size of the graph by a factor of two by looking only on the subgraph $\Gamma'$ induced by the neighbors of $v$.
Bob in his turn send Alice a name of the vertex $u$ in his independent set $I\cap \Gamma'$ which has degree at least $|V(\Gamma')|/2$.
Also in this case we can reduce the size of the remaining problem by a factor of two. Finally if both Alice and Bob can not send anything it is easy to see that 
$C\cap I = \emptyset$. By repeating this procedure at most $\log_2 m$ rounds, one can show that the deterministic
communication complexity satisfies $D(CL \text{-} IS_{\Gamma}) \leq O(\log_2^2 m)$. However, so far the best lower bound for this problem 
(see \cite{MR1723041}) is only asymptotically $2 \log_2 m$. 

For non-deterministic communication complexity of clique vs. independent set problem, it's easy to see that 
$N^{1}(CL \text{-} IS_{\Gamma})$  is always $\log m$. Indeed, for every vertex $v \in \Gamma$ consider 
the rectangle $R_v$ formed by all cliques vs. all independent sets containing $v$. By definition, these $m$ 
rectangles cover all $1$-inputs of the communication matrix $M$ of $CL \text{-} IS_{\Gamma}$.
On the other hand, determining the correct order of magnitude of $N^{0}(CL \text{-} IS_{\Gamma})$ 
is wide open except for the trivial lower bound $\log_2 m$. 
This lower bound follows from the simple fact that taking all single vertices as cliques vs. the same vertices as independent sets shows that 
the $m\times m$ identity matrix is a submatrix of $M$. 
Next we discuss the connection between the Alon-Saks-Seymour conjecture and
$CL\text{-}IS$ problem which was discovered by Alon and Haviv \cite{alonhaviv}.
This connection together with our counterexample gives a first nontrivial lower
bound for nondeterministic communication complexity of clique vs.
independent set problem. It implies that there exists a graph 
$\Gamma$ such that $N^{0}(CL \text{-} IS_{\Gamma}) \geq 6/5 \log_2 m -O(1)$.

Suppose we have a graph $G=(V,E)$, $V(G)=[n]$, $\bp(G)=m$ and a
partition of $E(G)$ as disjoint union of bicliques
$\{\mathcal{B}(U_i,W_i)\}_{i=1}^m$. Define the
characteristic vector $v_i$ of each biclique to be $v_i =
(v_{i1},\cdots, v_{in}) \in \{0,1,*\}^n$, so that
\begin{equation*}
v_{ij}=\begin{cases}
0 & \textrm{if}~j \in U_i \\
1 & \textrm{if}~j \in W_i\\
* & \textrm{otherwise}
\end{cases}
\end{equation*}

Using the notations above, we create a new graph $\Gamma$ on 
vertex set $[m]$. Two vertices $i$ and $i'$ are adjacent in $\Gamma$
if there exists $j \in [n]$ such that $v_{ij}=v_{i'j}=1$. Two
vertices $i$ and $i'$ are nonadjacent if there exists $j' \in [n]$
such that $v_{ij'}=v_{i'j'}=0$. In any other case,
arbitrarily assign an edge or non-edge between $i$ and $i'$. 
If there are two indices $j, j'$ such that
$v_{ij}=v_{i'j}=1$ and $v_{ij'}=v_{i'j'}=0$, then 
$j \in W_i \cap W_{i'}$ and 
$j' \in U_i \cap U_{i'}$. Therefore the edge $(j', j)$ is covered by two bicliques, which is impossible
since  $\cup_{i=1}^m \mathcal{B}(U_i,W_i)$ is an edge partition of $G$.
This shows that $\Gamma$ is well defined.

Now consider the $CL\text{-}IS$ problem on $\Gamma$. Define
$C_j=\{q\in [m]:v_{qj}=1\}$ and $I_j=\{q\in [m]:v_{qj}=0\}$. 
By definition of $\Gamma$, it is easy to see that $\{C_j\}$ are cliques and $\{I_j\}$ are
independent sets in this graph . Denote the
matrix of $CL \text{-} IS_{\Gamma}$ by $M$. Let $M'$ be a 
submatrix of $M$ corresponding to the rows determined by $\{C_j\}_{j=1}^n$ and
columns determined by $\{I_j\}_{j=1}^n$. Obviously $N^{0}(M) \geq
N^{0}(M')=\log_2 C^{0}(M')$. Assume that we have a covering of
$0$-entries of $M'$ by monochromatic rectangles, and let 
$R_1, \cdots, R_t$ be the rectangles
which cover the diagonal entries of $M'$. Note that if
$(p,q)$ is covered by $R_i$, then $M'_{pq}=M'_{qp}=0$ and thus $C_p
\cap I_q$ and $C_q \cap I_p$ are both empty. This implies that $(p,q)$ is not
an edge in graph $G$, since otherwise there must exist an index $i$ such that
$v_{ip}=0$, $v_{iq}=1$ or $v_{ip}=1$, $v_{iq}=0$. Then either $i \in I_p \cap C_q$ or
$i \in C_p \cap I_q$, which gives a contradiction.
In particular, the family of rectangles $\{R_i\}_{i=1}^t$
corresponds to a covering of graph $G$ by independent sets and
therefore $\chi(G) \leq t$. Thus we have that
$$N^{0}(M) \geq N^{0}(M')=\log_2 C^{0}(M') \geq \log_2 t \geq \log_2 \chi(G).$$ 
This estimate together with the existence of a graph $G$ (from Section 2) which has $\bp(G)=O(\chi(G)^{5/6})$,  
proves the following theorem.

\begin{thm}
There exists an infinite collection of graphs $\Gamma$, such that
$$N^{0}(CL\text{-} IS_{\Gamma}) \geq \frac{6}{5} \log_2 |V(\Gamma)| - O(1).$$
\end{thm}
In addition, the combination of the inequality $N^{0}(CL\text{-}IS_{\Gamma}) \geq \log_2 \chi(G)$ we just proved, and the result of  
Yannakakis  that $D(CL \text{-} IS_\Gamma) \leq O(\log_2^2 m)$, immediately gives a different derivation of the following result 
of Mubayi and Vishwanathan. It shows that if $\bp(G)=m$, then 
$$ \chi(G) \leq 2^{N^{0}(CL\text{-}IS_{\Gamma})} \leq 2^{D(CL \text{-} IS_\Gamma)} \leq 2^{O(\log_2^2 m)}.$$

From the above discussions, we know that any separation result between $\chi(G)$
and $\bp(G)$ gives corresponding separation between $N^{0}(CL \text{-} IS)$ and
the trivial lower bound $\log_2 |V(\Gamma)|$. 
We do not know whether the converse is also true yet. However, a weaker
converse does exist, as was observed by Alon and Haviv \cite{alonhaviv}. 
More precisely, the gap between
$N^{0}(CL \text{-} IS_{\Gamma})$ and $\log_2 |V(\Gamma)|$ implies a gap between
$\chi(H)$ and $2$-biclique partition number $\bp_2(H)$ for some graph $H$.

Let $\Gamma=(V,E)$ be a graph with vertices $V=\{v_1,\cdots, v_m\}$ and consider the following graph $H$.
The vertices of $H$ are all the pairs $(C,I)$ such that $C$ is a clique and $I$ is an independent
set in $\Gamma$ and $C \cap I = \emptyset$. Two vertices
$(C,I)$ and $(C',I')$ are adjacent if $C \cap I' \neq \emptyset$ or
$C' \cap I \neq \emptyset$. For every vertex $v_i$ in $\Gamma$, we define two subsets 
$U_i=\{(C,I): v_i \in C\}$ and $W_i=\{(C,I): v_i \in I\}$ of $H$.
These subsets have the following properties.

\begin{enumerate}
\item
$U_i$ and $W_i$ are disjoint.

\item
$(U_i,W_i)$ is a complete bipartite subgraph of $H$.

\item
$G'=\cup_{i=1}^m \mathcal{B}(U_i,W_i)$ and each edge of $H$ is
covered at most two times.
\end{enumerate}

The property $(1)$ holds since $C\cap I = \emptyset$ for any
vertex $(C,I)$ of $H$. To verify $(2)$, consider two vertices $(C,I) \in
U_i$ and $(C',I') \in W_i$. Then $v_i \in C \cap I'$, which means $C \cap
I' \neq \emptyset$ and thus $(C,I)$ and $(C',I')$ are adjacent in
$H$. To prove $(3)$, note that by definition, any edge $(C,I) \sim (C',I')$ in $G'$ either satisfies
$C \cap I'\neq \emptyset$ or $C' \cap I \neq \emptyset$ or both. If $C \cap
I'\neq \emptyset$, then there is a unique $i$ (since $|C \cap I'| \leq 1$)
such that $v_i \in C$ and $v_i \in I'$, which means that this edge belongs
to $\mathcal{B}(U_i,W_i)$. The similar conclusion holds in the case when $C' \cap
I\neq \emptyset$. Thus every edge of $H$ is covered by $\{\mathcal{B}(U_i,W_i)\}_{i=1}^m$ 
either once or twice. This shows that $\bp_2(H) \leq m=|V(\Gamma)|$. 

Next we bound the chromatic number of $H$ from below by a function of
$N^{0}(CL \text{-} IS_{\Gamma})$. Denote the matrix of $CL\text{-}IS_{\Gamma}$ by $M$.
By definition, an independent set $I'=\{(C_1,I_1), \cdots, (C_l,I_l)\}$ of graph $H$ 
corresponds to an all-zero submatrix of $M$, whose rows and columns are indexed by
$C_1, \cdots, C_l$ and $I_1, \cdots, I_l$ respectively. Thus a proper
coloring of $H$ corresponds to a covering of the $0$-entries
of $M$ by monochromatic rectangles. Therefore $\chi(H) \geq C^{0}(M)=C^{0}(CL \text{-} IS_{\Gamma}) \geq 2^{N^{0}(CL \text{-}
IS_{\Gamma})}$ and hence we established the following claim.

\begin{claim}
For every graph $\Gamma$ there exists a graph $H$ such that
$$\bp_2(H)\leq |V(\Gamma)| ~~~\mbox{and}~~~ \chi(H) \geq 2^{N^{0}(CL \text{-} IS_{\Gamma})}.$$
\end{claim}

\section{Concluding remarks}
\label{section_concluding}

In this paper we constructed a graph which has a polynomial gap between the chromatic number and
the biclique partition number, thereby disproving the Alon-Saks-Seymour conjecture. A very interesting problem which remains widely open
is to determine how large this gap can be. In communication complexity
it is a long standing open problem to prove an $\Omega (\log^2 N)$ lower bound on the complexity of clique vs. independent
set problem for graph on $N$ vertices. Since, as we already explained in the previous section, this problem is closely related to
the Alon-Saks-Seymour conjecture, it is plausible to believe that one can obtain a corresponding gap between
chromatic and biclique partition numbers. We conjecture that there exists a graph $G$ with biclique partition number $k$ and chromatic number 
at least $2^{c \log^2 k}$, for some constant $c>0$. Existence of such graph will also resolve the complexity of clique vs. independent set problem.

Another intriguing question which deserves further study is to determine the $t$-biclique covering numbers of complete graphs.
This will also solve the problem of the maximum possible cardinality of $t$-neighborly family of standard boxes in finite dimensional
Euclidean spaces. Even the asymptotics of $\bp_t(K_k)$ is only known up to a multiplicative constant factor. In the first open case when $t=2$,
the best current bounds are $ (1+o(1)) k^{1/2} \leq \bp_2(K_k) \leq (1+o(1)) 2k^{1/2}$ and it would be interesting to
close this gap.

\vspace{0.25cm}
\noindent
\textbf{Acknowledgment.} We would like to thank Noga Alon
for explaining to us his results with I. Haviv  on connection between the Alon-Saks-Seymour
conjecture and the clique vs. independent set problem.

\end{document}